\documentclass[a4paper,11pt]{article}
\usepackage[russian,english]{babel}
\usepackage{amsfonts,amsmath,amssymb,theorem}
\usepackage{hyperref} 
\usepackage{colortbl}
\usepackage{ifthen}
\usepackage{hhline}

\makeatletter
\def\@seccntformat#1{\csname the#1\endcsname.\ } 
\def\@biblabel#1{#1.} 
\makeatother

\date{}

    \newif\ifNoRemark
    \def\addtheorem#1#2#3#4{
    \ifthenelse{\expandafter\isundefined\csname the#2\endcsname}{\newcounter{#2}}{}
    \newenvironment{#1}[1][\global\NoRemarktrue]
     {\par\addvspace{2mm plus 0.5mm minus 0.2mm}\noindent
       \refstepcounter{#2}{\bf #3~\csname the#2\endcsname
      \vphantom{##1}\ifNoRemark.\ \else\ (##1).\fi}\begingroup #4}%
      {\endgroup\par\addvspace{1mm plus 0.5mm minus 0.2mm}\global\NoRemarkfalse}
    \expandafter\newcommand\csname b#1\endcsname{\begin{#1}}
    \expandafter\newcommand\csname e#1\endcsname{\end{#1}}
    }

\addtheorem{theorem}{thrm}{Theorem}{\sl}
\addtheorem{lemma}{lmm}{Lemma}{\sl}
\addtheorem{pro}{prpstn}{Proposition}{\sl}
\addtheorem{corol}{corol}{Corollary}{\sl}

\begin{document}

\title{On the connection between correlation-immune functions and
 perfect $2$-colorings of the Boolean  $n$-cube  \thanks{The work is supported by RFBR (grants
10-01-00424, 10-01-00616)}}
\author{Vladimir N.~Potapov}

\maketitle

\begin{center}
\textit{Sobolev Institute of Mathematics,\\
Novosibirsk State University, Novosibirsk }
\end{center}

\begin{abstract}

A coloring of the Boolean $n$-cube  is called perfect if, for every
vertex $x$,  the collection of the colors of the neighbors of $x$ depends only on
the color of $x$.  A Boolean function is called correlation-immune of
degree $n-m$ if it takes the value $1$ the same number of  times for each
$m$-face of the Boolean $n$-cube. In the present paper it is proven that
each Boolean function $\chi^S$ ($S\subset E^n$) satisfies the
inequality
  $${\rm nei}(S)+ 2({\rm cor}(S)+1)(1-\rho(S))\leq
n,$$ where ${\rm cor}(S)$
  is the maximum degree of the correlation immunity of $\chi^S$, ${\rm nei} (S)= \frac{1}{|S|}\sum\limits_{x\in
S}|B(x)\cap S|-1$ is the average number of neighbors in the set $S$
for vertices in $S$, and
  $\rho(S)=|S|/2^n$ is the density of the set $S$.
  Moreover, the function $\chi^S$ is a perfect coloring if and only
  if we obtain an equality in the above formula.


\end{abstract}


Let $E^n=\{ 0,1\}^n$ be the $n$-dimensional Boolean cube ($n$-cube).
  Define the operation   $[x,y]=(x_1y_1,\dots,x_ny_n)$ and
the inner product $\langle x,y\rangle=x_1y_1\oplus\dots\oplus x_ny_n$
for vectors $x,y\in E^n$. The number of unities in a vector  $y\in
E^n$ is called the {\it weight} of $y$  and is denoted by  $wt(y)=
\langle y, \overline{1}\rangle$. The sets $E^n_y(z)=\{x\in E^n :
  [x,y]=[z,y]\}$, where $wt(y)=m$, are called $(n-m)$-{\it faces}.

Let $S\subset E^n$  and let $\chi^S$ be the characteristic function of $S$.
The function $\chi^S$ is said to be {\it correlation-immune
of order} $n-m$ if, for all $m$-faces $E^n_y(z)$, the intersections
 $E^n_y(z)\cap S$ have the same cardinality.
  The maximum order of the correlation immunity of the function $\chi^S$
  is denoted by ${\rm cor}(S)$,   ${\rm cor}(S)=\max \{n-m\}$. The number
  $\rho(S)=|S|/2^n$ is the {\it
density} of the function $\chi^S$. We will always assume that
$\rho(S)\leq 1/2$, since otherwise we may consider
$E^n\setminus S$ instead of $S$. If $\rho(S)=1/2$ then the
correlation immune  function $\chi^S$ is called  {\it balanced}.

The {\it Hamming distance} $d({x},{y})$ between two vectors ${x},
{y}\in \Sigma^n$ is the number of positions at which they differ. The
unit ball $\{y\in E^n : d(x,y)\leq 1\}$ is denoted by ${B}({x})$.
Define the {\it neighbors number} ${\rm nei} (S)$ to be the average number of
neighbors in the set $S$ for vertices in $S$, i.\,e., ${\rm nei}
(S)= \frac{1}{|S|}\sum\limits_{x\in S}|B(x)\cap S|-1$.

A mapping $Col:E^n\rightarrow\{1,\dots,k\}$ is called a {\it perfect
 coloring} with matrix of parameters $A=\{a_{ij}\}$  if, for all $i$,
  $j$, for every  vertex of color $i$, the number of  its neighbors of color $j$ is
  equal to $a_{ij}$.
In what follows we will only consider colorings in two colors. Moreover,
for convenience we will assume that the set of colors is $\{0,1\}$.
In this case the function $Col$ is Boolean and
$Col=\chi^S$, where $S$ is the set of 1-colored vertices.

   A {\it perfect code} $C\subset E^n$ can be defined as a perfect coloring with
   matrix of parameters $  A= \left(
 \begin{array}{ccc}
0 & n  \\
  1  & n-1  \\
  \end{array} \right).
$ Colorings with such parameters exist only if   $n=2^m-1$ ($m$ is
integer). A list of accessible parameters and corresponding
constructions of  perfect colorings can be found in \cite{Fl1} and
\cite{Fl3}.

It is well known (see \cite{Fl2}) that a perfect coloring of
the $n$-cube with matrix of parameters $\left(
 \begin{array}{ccc}
n-b & b  \\
  c  & n-c  \\
  \end{array} \right)
$ is  a correlation-immune function of degree $\frac{b+c}{2} -1$.
Therefore, if the vertices of some set $S$ are regularly distributed on balls
then the vertices of the set  are uniformly distributed on faces.
It is of interest to clarify the possibility of the reverse implication.

  In \cite{Fl2} it is established  that for each unbalanced Boolean
  function $\chi^S$ ($S\subset E^n$)
  the inequality
  ${\rm cor}(S)\leq \frac{2n}{3} -1$  holds. Moreover, in the case of the equality
    ${\rm cor}(S)= \frac{2n}{3} -1$, the function $\chi^S$  is a perfect coloring.
Similarly, if for any set $S\subset E^n$ the Bierbrauer-Friedman
inequality \cite{Fr,Br} $\rho(S)\geq 1-\frac{n}{2({\rm cor}(S)+1)}$
becomes an equality then the function $\chi^S$ is a perfect
$2$-coloring (see \cite{Pot}). Consequently, in the extremal cases, the regular
distribution on balls follows from the uniform distribution on
faces. In the present paper we prove the following theorem:

\btheorem\label{thki1}

 {\rm (a)} For each Boolean function
   $a=\chi^S$, where $S\subset E^n$ and $\rho(S)\leq
1/2$, the inequality ${\rm nei}(S)+ 2({\rm cor}(S)+1)(1-\rho(S))\leq
n$   holds.

 {\rm (b)} A Boolean function $a=\chi^S$ is  a perfect $2$-coloring
 if and only if\\
${\rm nei}(S)+ 2({\rm cor}(S)+1)(1-\rho(S))=n$. \etheorem

In the proof of the theorem we employ the idea from the papers
\cite{Fr} and \cite{OP}. Before proceeding to the proof we introduce
some necessary concepts.

The set $\mathbb{V}$ of all functions $a:E^n\rightarrow\mathbb{Q}$ has the natural
structure of an $2^n$-dimensional vector space. It is well
known that the functions $f^v(u)=(-1)^{\langle u,v\rangle}$, where
$v\in E^n$, constitute a Fourier  orthogonal basis of the space $\mathbb{V}$.
The  {\it Fourier transform} $\widehat{a}$ of
a function $a$ is determined by the equality
$$\widehat{a}(v)= \sum\limits_{u\in E^n}a(u)(-1)^{\langle
u,v\rangle}.
$$
 Obviously, the function $\widehat{a}$ is the inner product
of the vectors $a$ and $f^v$ in $\mathbb{V}$. Since $\langle f^v,
f^v\rangle =2^n$ for every vertex $v\in E^n$, we have the equalities
\begin{equation}\label{ki5}
a(u)= \frac{1}{2^n}\sum\limits_{v\in E^n}\widehat{a}(v)(-1)^{\langle
u,v\rangle}, \ \ \ \
\end{equation}
\begin{equation}\label{ki6}
 2^{n}\sum\limits_{u\in E^n}a^2(u)= \sum\limits_{v\in
E^n}\widehat{a}^2(v).
\end{equation}

We need the following well-known statements.

\bpro\label{stki1} (\cite{Fl2}, \cite{MW}) A Boolean function
$a=\chi^S$ is correlation immune of degree  $m$ if and only if
$\widehat{a}(v)=0$ for every $v\in E^n$ such that $0<wt(v)\leq
m$.\epro

\bpro\label{stki2} (\cite{Fl1})

{\rm (a)} Let $a=\chi^S$ be a perfect coloring with matrix of
parameters $\left(
 \begin{array}{ccc}
n-b & b  \\
  c  & n-c  \\
  \end{array} \right).
$ Then $\widehat{a}(v)=0$ for every $v\in E^n$ such that $wt(v)\neq
0,\frac{b+c}{2}$.

{\rm (b)} Let $a=\chi^S$ be a Boolean function. If
$\widehat{a}(v)=0$ for every vertex $v\in E^n$ such that $wt(v)\neq
0,k$ then $a$ is a perfect coloring. \epro

\bcorol\label{corki0} Let $a=\chi^S$ be a perfect coloring with
matrix of parameters $\left(
 \begin{array}{ccc}
n-b & b  \\
  c  & n-c  \\
  \end{array} \right).
$ Then
 ${\rm cor}(S)= \frac{b+c}{2} -1$.\ecorol

Put $a'(v)=\widehat{a}(v)/|S|$. From  (\ref{ki6}) we obtain
\begin{equation}\label{ki7}
\sum\limits_{v\in E^n} (a'(v))^2=
\frac{2^{n}}{|S|^2}\sum\limits_{u\in E^n}a^2(u)= 2^n/|S|.
\end{equation}

Given a Boolean function $a=\chi^S$,  the {\it weight distribution }
$(B_0(a),\dots,B_n(a))$ is determined by  $$B_i(a)=\frac{1}{|S|}
|\{v,u\in S \ |\ wt(u+v)=i\}|.$$ Obviously, $B_1(a)={\rm nei}(S)$,
where $a=\chi^S$.

The {\it MacWilliams transform} of a weight distribution
$(B_0(a),\dots,B_n(a))$ is the array\\ $(B'_0(a),\dots,B'_n(a))$,
where $B'_k(a)=\frac{1}{|S|}\sum_{i=0}^nB_i(a)P_k(i)$ and $P_k$ are
the Kravchuk polynomials.

In \cite{MW} the following equalities are proven:
\begin{equation}\label{ki10}
B_k(a)= \frac{|S|}{2^n}\sum_{i=0}^nB'_i(a)P_k(i),
\end{equation}
\begin{equation}\label{ki11}
B'_k(a)=\sum_{wt(v)=k}(a'(v))^2.
\end{equation}
From (\ref{ki7}) and (\ref{ki11}) we obtain the following:

\bcorol\label{corolki1} If $\overline{0}\in S$ and $a=\chi^S$ then

{\rm (a)} $B'_k(a)\geq 0$ for $i=0,\dots,n$;

{\rm (b)} $B'_k(a)= 0$ $\Leftrightarrow$ $a'(v)=0$ for every vector
$v\in E^n$ with weight $wt(v)=k$;

{\rm (c)} $B'_0(a)=1$;

{\rm (d)} $\sum\limits_{i=0}^nB'_k(a)=2^n/|S|$.

\ecorol

Propositions \ref{stki1}, \ref{stki2} and Corollary
\ref{corolki1} (b) imply the following well-known statements
(see \cite{Fl2} and \cite{OP}).

\bcorol\label{corolki2} If $\overline{0}\in S$ and $a=\chi^S$ then
$B'_k(a)=0$ for $0<k\leq {\rm cor}(S)$.\ecorol

\bcorol\label{corolki3} Let $\overline{0}\in S$, $a=\chi^S$, and
let $B'_i(a)=0$  for $i\neq 0,k$. Then  $\chi^S$ is a perfect
coloring.\ecorol

{\bf Proof of the theorem.} Without lost of generality we suppose
that $\overline{0}\in S$. Put $t={\rm cor}(S)$.  From (\ref{ki10}),
Corollary \ref{corolki1}(c), and Corollary \ref{corolki2}
we obtain the equality
$${\rm nei}(S)/\rho(S)= P_1(0)+\sum\limits_{i=t+1}^nB'_i(a)P_1(i).$$
Since  $P_1(i)=n-2i$ (see \cite{OP} or \cite{MW}), we have
$${\rm nei}(S)/\rho(S)\leq n + \sum\limits_{i=t+1}^nB'_i(a)(n-2i).$$
Corollary \ref{corolki1} (a),(d) implies
$\sum\limits_{i=0}^nB'_k(a)=1/\rho(S)$ and $B'_k(a)\geq 0$. Hence,
$${\rm nei}(S)\leq \rho(S)n + (n-2(t+1))\left(1-{\rho(S)}\right).$$
Moreover, the equality
\begin{equation}\label{ki17}
{\rm nei}(S)= \rho(S)n + (n-2(t+1))\left(1-{\rho(S)}\right)
\end{equation} holds if and only if
$B'_i(a)=0$ for $i\geq t+2$. Then from Corollary \ref{corolki3}
(b) we conclude that $\chi^S$ is a perfect coloring.

Each perfect $2$-coloring satisfies (\ref{ki17}),
which is a consequence of Proposition \ref{stki2} (b) and  Corollary
\ref{corolki3} (b).
 $\bigtriangleup$

For perfect codes,  a similar  theorem was previously proven in
\cite{OP}. Namely, if  ${\rm cor}(S)={\rm cor}(H)$ and
  $\rho(S)=\rho(H)$, where
  $S,H\subset E^n$ and  $H$  is a perfect code, then $S$ is also a perfect code.

\begin{center}

\end{center}

\end{document}